\documentstyle[12pt]{amsart}
\input{amssymb.sty}
\input xy
\xyoption{all}
\theoremstyle{plain}
  \newtheorem{thm}{Theorem}[section]
  \newtheorem{lem}[thm]{Lemma}
  \newtheorem{cor}[thm]{Corollary}
  \newtheorem{prop}[thm]{Proposition}
\theoremstyle{definition}
  \newtheorem{defn}[thm]{Definition}

  \newtheorem{clm}[thm]{Claim}
  \newtheorem{notation}{Notation\!\!}

\theoremstyle{remark}
  
\newcommand{\Coh}{{\operatorname{Coh}}}

\newcommand{\Ext}{\operatorname{Ext}}

\newcommand{\Hom}{\operatorname{Hom}}
\newcommand{\id}{\operatorname{id}}
\newcommand{\IIm}{\operatorname{Im}}
\newcommand{\Ker}{\operatorname{Ker}}

\newcommand{\Mod}{\operatorname{Mod}}

\newcommand{\NS}{\operatorname{NS}}

\newcommand{\QQ}{{\mathbb{Q}}}

\newcommand{\rk}{{\operatorname{rk}}}

\newcommand{\SL}{\operatorname{SL}}
\newcommand{\Spec}{\operatorname{Spec}}

\newcommand{\ZZ}{{\mathbb{Z}}}

\newcommand{\bA}{\overline{A}}
\newcommand{\bE}{\overline{E}}
\newcommand{\bEE}{\overline{{\cal E}}}
\newcommand{\bGam}{\overline{\Gamma}}
\newcommand{\bT}{\overline{T}}
\newcommand{\bV}{\overline{V}}
\newcommand{\cc}{{\frak c}}
\newcommand{\dotimes}{\underline{\otimes}}
\newcommand{\EE}{{\cal E}}
\newcommand{\FFl}{{\cal F}^{(l)}}
\newcommand{\FFr}{{\cal F}^{(r)}}
\newcommand{\Fl}{F^{(l)}}
\newcommand{\Fr}{F^{(r)}}
\newcommand{\Hf}{{\frak f}}
\newcommand{\HomSF}{\Hom_{SF}}
\newcommand{\FF}{{\cal F}}
\newcommand{\Gamdot}{\Gamma_{\bullet}}
\newcommand{\Gaml}{\Gamma^{(l)}}
\newcommand{\Gamr}{\Gamma^{(r)}}
\newcommand{\GG}{{\cal G}}

\newcommand{\Ibr}{I^{\bullet(r)}}
\newcommand{\Ms}{M^s}
\newcommand{\Mm}{M_-}
\newcommand{\Mp}{M_+}
\newcommand{\ob}{\operatorname{ob}}
\newcommand{\tilM}{\tilde{M}}
\newcommand{\Vl}{V^{(l)}}
\newcommand{\Vr}{V^{(r)}}
\begin{document}
\bibliographystyle{amsplain}
\title[blowing-ups describing polarization change]
{Blowing-ups describing the polarization change of
moduli schemes of semistable sheaves \\
of general rank}
\author{Kimiko Yamada}
\address{Department of Mathematics, Faculty of science and technology,
Sophia Universiy, 7-1, Kioi-cho, Chiyoda-ku, Tokyo, Japan}
\email{yamada@@mm.sophia.ac.jp}
\thanks{The author was supported by Japan Society for the Promotion of Science.}

\maketitle

\begin{abstract}
Let $H$ and $H'$ be two ample line bundles over a smooth projective
surface $X$, and $M(H)$ (resp. $M(H')$) the coarse moduli scheme of
$H$-semistable (resp. $H'$-semistable) sheaves of fixed type
$(r,c_1,c_2)$.
We construct a sequence of blowing-ups which describes how
$M(H)$ differs from $M(H')$ not only when $r=2$ but also when $r$ is arbitrary.
Means we here utilize are elementary transforms and the notion of a sheaf with flag.
\end{abstract}

\section{Introduction}

Let $X$ be a nonsingular projective surface over an algebraically
closed field $k$ with character zero, $H_-$ and $H_+$ two ample line
bundles over $X$, and $\cc=(r,c_1,c_2)$ an element of $\ZZ\times
\NS(X)\times\ZZ$.
There exists the coarse moduli scheme $M(H_-,\cc)$, which is
projective over $k$, of S-equivalence classes of $H_-$-semistable
sheaves $E$ on $X$ such that $(r(E),c_1(E),c_2(E))=\cc$ by
\cite{Gi:moduli}.
When $2rc_2-(r-1)c_1^2$ is sufficiently large, $M(H_-,\cc)$
and $M(H_+,\cc)$ are birationally equivalent from
\cite[Theorem 4.C.7]{HL:text}. 
With this in mind, we shall construct a sequence of morphisms
\begin{equation}\label{eq:MainDiag}
\xymatrix{
 & \tilM_0 \ar[dl]_{p_0} \ar[dr]^{q_0} & &  \dots \ar[dl]_{p_1} &
 \tilM_{m-1} \ar[dr] ^{q_{m-1}} & \\
 M(\Delta_0,\cc) \ar[d]^h & & M(\Delta_1,\cc)  &  & &
 M(\Delta_m,\cc) \ar[d]^h  \\
 M(H_-,\cc) &  & & & & M(H_+,\cc)}
\end{equation}
assuming that $H_-$ and $H_+$ lie in adjacent chambers
(\cite[Definition 2.1]{Yos:chamber}) of type $\cc$ . \par
To execute our purpose we utilize elementary transforms and introduce
a sheaf with flag, or a SF for short.
Elementary transforms have appeared in the study of the polarization
change problem for stability conditions of rank-two stable sheaves.
However we can not directly apply this way to the general-rank case
partly because an $H_-$-semistable and not $H_+$-semistable sheaf
of type $\cc$ is $H_-$-stable if its rank is two, but it is not
necessarily $H_-$-stable in general. For example, if a sheaf
$F$ of type $\cc$ is $H_-$-semistable and not $H_+$-semistable,
then $F\oplus F$ is $H_-$-semistable, not $H_+$-semistable and
not $H_-$-stable.
It is unfavorable since the complement of
$M^s(H_-,\cc)\subset M(H_-,\cc)$, the open set of all
$H_-$-stable sheaves, is complicated.
Hence in Section \ref{ss:SF} we introduce a sheaf with flag $($SF$)$ 
and
its $\Delta$-stability with respect to $(L,C)$, where $\Delta$ is
a parameter, $L$ a line bundle on $X$ and $C\subset X$ an effective
divisor.
As is discussed in Section \ref{ss:moduli} the coarse moduli scheme of
$\Delta$-semistable SFs exists; $M(\Delta,\cc)$ at
\eqref{eq:MainDiag} is deduced from it.
Corollary \ref{cor:FlFr} shows that under some condition 
the problem of observing how
stability conditions of SFs vary as parameters $\Delta$ do is
similar to the polarization change problem for rank-two
stable sheaves.
With this corollary as a base, we obtain blowing-ups $p_i$ and $q_i$ at
\eqref{eq:MainDiag}, whose centers are topics in Section \ref{ss:SetP}.
The morphism $h:M(\Delta_0,\cc)\rightarrow M(H_-,\cc)$
at \eqref{eq:MainDiag} is naturally induced when $\Delta_0$ and $(L,C)$
are chosen appropriately.
Its restriction 
$h: h^{-1}(M^s(H_-,\cc))\rightarrow M^s(H_-,\cc)$ is
a Grassmannian-bundle in \'{e}tale topology. \par
Here let us mention the background.
With its relation to the wall-crossing formula of Donaldson polynomials,
the polarization change problem for stability conditions of sheaves has
been a subject of study.
Matsuki-Wentworth \cite{MW:Mumford} pointed out in general-rank case
that this problem is a subject concerning how the GIT quotient of a
quasi-projective scheme $S$ by a reductive group $G$ varies as 
$G$-linearized ample line bundles of $S$ do, and connected
$M(H_-)$ and $M(H_+)$ by a sequence of Thaddeus-type
flips (\cite{Tha:GITflip}).
On the other hand,
elementary transforms, refer to \cite{Mar:ElemTrans}
and \cite[Appendix]{Fr:elliptic} about general information,
has the following advantages:
(i) birational transforms obtained there are blowing-ups whose centers
are derived by a canonical, moduli-theoretic way;
(ii) when two parameters $\alpha$ and $\alpha'$
defining stability conditions of objects are given, one not only
connects the moduli scheme of $\alpha$-semistable objects with that
of $\alpha'$-semistable ones, but also relates their universal
families, if exist.
Ellingsrud-G\"{o}ttsche \cite{EG:variation} and Friedman-Qin
\cite{FQ:flips} proposed to apply elementary transform to the case
where $r=2$ and the wall of type $\cc$ separating $H_-$ and $H_+$
is good, so the natural subset
\[ M(H_-,\cc)\supset P=\left\{[E] \bigm| 
 \text{$E$ is not $H_+$-semistable} \right\} \]
is relatively easy to handle.
These papers have stimulated the author to write this article.
The author aims to consider this problem with no restriction on this
wall. As a result this subset $P$ unkindly behaves in general,
and we have to observe its (infinitesimal) structure in more detail.
The preceding paper \cite{Yam:Dthesis} dealt with rank-two case and
this article the case where $r$ is arbitrary, and hence
we need further devices explained above.
We finally note that while writing this article the author found that
also Mochizuki used the notion of sheaves with flag in
\cite{Mo:invariants}, where he considered not birational transforms
describing the variation of moduli schemes, but the wall-crossing
formula of Donaldson polynomials in general-rank case. \par
The content of this article is as follows.
In Section \ref{ss:SF} we define some basic terms and show Corollary
\ref{cor:FlFr} mentioned above.
In Section \ref{ss:moduli} we construct the moduli scheme
$M(\Delta,\Hf)$ of $\Delta$-semistable SFs of type $\Hf$
and study its infinitesimal structure.
The scheme $M(\Delta,\cc)$ at \eqref{eq:MainDiag} is 
the union of some connected
components of $M(\Delta,\Hf)$.
We focus in Section \ref{ss:SetP} on the subscheme
$P\subset M(\Delta_-,\Hf)$ consisting of SFs which are
$\Delta_-$-semistable and not $\Delta_+$-semistable, and discuss its
relative obstruction theory.
In Section \ref{ss:blowup} we arrive at the sequence \eqref{eq:MainDiag}
by elementary transforms.
\begin{notation}
$X$ is a smooth projective surface over an algebraically closed field
$k$ with character zero.
For $k$-schemes $S$ and $T$, $T_S$ means $T\times S$ and
$p_T: T_S \rightarrow S$ is the natural projection.
${\cal H}^i(A)$ is the $i$-th cohomology of $A\in D(T):=D(Qcoh(T))$.
$Mc(A\rightarrow B)$ is the mapping cone of a morphism $A\rightarrow B$
in $D(T)$. $\dotimes$ stands for the derived functor of $\otimes$.
\end{notation}
\section{Sheaves with flag }
\label{ss:SF}
\begin{defn}\label{defn:SF}
{\it A sheaf with flag (SF) of length} $n$ is a pair
$\EE=(E, \left\{\Gamma_i \right\}_{i=1}^n)$ consisting of a coherent
sheaf $E$ on $X$ and a flag of vector spaces
$\Gamma_1 \subset \Gamma_2 \subset \dots \Gamma_n \subset H^0(E)$.
%{\it The type of a SF} $\EE=(E,\Gamma_{\bullet})$ is the
%class 
%\begin{equation*}
%\left( r(E), c_1(E), c_2(E), \rk\Gamma_1, \dots, \rk\Gamma_n \right)
%\in \NN\times\NS(X)\times\ZZ\times\ZZ^{\times n}.
%\end{equation*}
%
{\it A homomorphism} $f: \EE'=\left(E',\Gamma'_{\bullet}\right)\rightarrow
\EE=(E,\Gamma_{\bullet})$ {\it of SFs of length} $n$ is a homomorphism
$f: E'\rightarrow E$ of sheaves which preserves their flag structures.
$\HomSF(\EE',\EE)$ denotes the set of all
homomorphisms $f:\EE'\rightarrow \EE$ of SFs.
We say that a sequence
$\EE^{(0)} \overset{f^{(0)}}{\rightarrow} \EE^{(1)} 
 \overset{f^{(1)}}{\rightarrow} \EE^{(2)}$ of SFs
$\EE^{(j)}=(E^{(j)},\Gamdot^{(j)})$ and homomorphisms is {\it exact}
if both $E^{(0)} \rightarrow E^{(1)} \rightarrow E^{(2)}$ and
$\Gamma_i^{(0)} \rightarrow \Gamma_i^{(1)} \rightarrow \Gamma_i^{(2)}$
($i$ is arbitrary) are exact.
{\it A sub SF} $\EE'\subset\EE$ is given by a homomorphism
$\iota:\EE'\rightarrow \EE$ of SFs such that $\iota :E' \rightarrow E$ is
injective.
A sub SF
 $\EE'=\left(E',\Gamma'_{\bullet}\right)\subset\EE=(E,\Gamma_{\bullet})$
is said to be {\it saturated} if the induced 
homomorphism
$\Gamma_i / \Gamma'_i\rightarrow H^0(E/E')$ is injective for all $i$;
in other words, $\Gamma'_i=H^0(E')\cap\Gamma_i$ for all $i$.
In this case, also $\EE/\EE'=\left(E/E', \left\{\Gamma_i/\Gamma'_i\right\}_i
\right)$ is a SF.
A SF $\EE=\left(E, \Gamma_{\bullet}\right)$ of length $n$ is {\it full}
if it holds that $\rk\Gamma_i=i$ for all $i$ and that $n=h^0(E)$.
\end{defn}
\begin{defn}
Let ${\cal O}(1)$ be an ample line bundle on $X$, $L$ a line bundle on
$X$, and $C\subset X$ an effective divisor on $X$.
A sheaf $E$ is said to be {\it of type} $\Hf'\in \QQ[l]^{\times 3}$ if 
\[  \bigl( \chi(E(l)),\, \chi(E\otimes L(-C)(l)),\, \chi(E\otimes L(l))
\bigr) =\Hf',\]
and a SF $\EE$ of length $n$ is said to be {\it of type} 
$\Hf \in \QQ[l]^{\times 3}\times \ZZ^{\times n}$ if
\[ \bigl( \chi(E(l)),\, \chi(E\otimes L(-C)(l)),\, \chi(E\otimes L(l)),\, 
\rk\Gamma_1, \dots, \rk\Gamma_n \bigr) =\Hf.\]
When a parameter
$\Delta=(a, \delta_1, \dots, \delta_n)\in (0,1)\times\QQ^{\times n}_{>0}$
is given,
we also define the {\it reduced Hilbert polynomial} of a SF
$\EE$ of length $n$ with $\rk E>0$ by
\[ p^{\Delta}(\EE)(l)=\frac{1}{\rk E}\left\{(1-a)\chi
\left(E\otimes L(-C)(l)\right)+a\chi\left(E\otimes L(l)\right)
 + \sum_{i=1}^{n} \delta_i \cdot\rk\Gamma_i \right\}\in\QQ[l]. \]
\end{defn}
\begin{defn}
For a parameter
$\Delta=(a, \delta_1, \dots, \delta_n)\in (0,1) \times\QQ^{\times n}_{>0}$, 
we say that a SF $\EE=(E, \Gamma_{\bullet})$ of length $n$ is 
{\it $\Delta$-stable $($resp. semistable$)$} if $E$ is torsion-free and
it holds that
$ p^{\Delta}(\EE')<p^{\Delta}(\EE) $ (resp. $\leq$) for any proper sub SF 
$\EE'\subset\EE$. 
We define the {\it S-equivalence} of
$\Delta$-semistable SFs in the same way as the case of
semistable sheaves \cite[p. 22]{HL:text}.
\end{defn}
For $\Hf\in Q[l]^{\times 3}\times \ZZ^{\times n}$,
we set ${\cal S}_1(\Hf)$ to be the set of all nonzero SFs $\EE'$ of length
$n$ such that there are a SF $\EE=(E,\Gamdot)$ of type $\Hf$ and 
a parameter $\Delta_0$ satisfying
(i) $E$ is ${\cal O}(1)$-semistable, (ii) $\EE'$ is a proper sub SF of
$\EE$ and (iii) $p^{\Delta_0}(\EE')=p^{\Delta_0}(\EE)$.
Any SF $\EE'\in{\cal S}_1(\Hf)$ gives a subset in
$(0,1)\times \QQ^{\times n}_{>0}$
\[W(\EE',\Hf) =\left\{ \Delta=(a,\delta_{\bullet})
\bigm| \text{ $p^{\Delta}(\EE')=p^{\Delta}(\EE)$
for any SF $\EE$ of type $\Hf$ } \right\}.  \]
%One can check this is a hyperplane.
Grothendieck's lemma on boundedness \cite[p. 29]{HL:text} implies
$\left\{ W(\EE',\Hf) \bigm| \EE'\in {\cal S}_1(\Hf) \right\}$ is finite.
\begin{defn}
This $W(\EE',\Hf)$ is called a {\it SF-wall of type}
$\Hf$ if it is a proper subset of $(0,1)\times \QQ^{\times n}_{>0}$.
{\it A SF-chamber of type} $\Hf$ is a connected component of
the complement of the union of all SF-walls of type $\Hf$.
$\Delta$-semistability of SFs of type $\Hf$ does not change
unless $\Delta$ passes through a SF-wall of type $\Hf$.
\end{defn}
We say that a ${\cal O}_X$-module $F$ has {\it the property} $(O)$
(resp. $(O_m)$) with respect to an ample line bundle ${\cal O}(1)$
if $F$ (resp. $F(m)$) is generated by global sections and its higher 
cohomologies vanish.
For $\Hf'\in \QQ[l]^{\times 3}$
we define two families as follows:
let ${\cal S}_2(\Hf') $ be the set of all
${\cal O}(1)$-slope-semistable sheaves of type $\Hf'$ on $X$, and
let ${\cal S}_3(\Hf')$ be the set of all sheaves $E'$ on $X$ such that
$E'$ is a subsheaf of a certain $E\in{\cal S}_2(\Hf')$ with torsion-free
quotient and satisfies $\mu_{{\cal O}(1)}(E')=\mu_{{\cal O}(1)}(E)$.
If one replace $E$ with $E(m)$ where $m$ is sufficiently large, 
then it holds that
\begin{equation}\label{eq:(O)}
\text{Every member of ${\cal S}_2(\Hf')\cup {\cal S}_3(\Hf')$ has
the property $(O)$.}
\end{equation}
\begin{defn}\label{defn:(A)}
We say that
$\Hf=(\Hf', l_1,\dots, l_n)\in \QQ[l]^{\times 3}\times \ZZ^{\times n}$
has {\it the property} $(A)$ if
\eqref{eq:(O)} is valid for $\Hf'=(f,f_0,f_1)$ and if
$n$ and $\Hf$, respectively, equal $f(0)$ and $(\Hf',1,2,\dots, n)$.
\end{defn}
One can verify that if $\Hf$ has the property $(A)$ and
if a parameter $\Delta$ is contained in no SF-wall of type $\Hf$,
then a $\Delta$-semistable SF $\EE$ of type $\Hf$ is always $\Delta$-stable.
Moreover, we have the proposition below.
\begin{prop}\label{prop:Rk2-like}
Assume that $\Hf$ has the property $(A)$.
Suppose that a parameter $\Delta_0$ is contained in just one SF-wall
of type $\Hf$ and a SF $\EE$ of type $\Hf$ is $\Delta_0$-semistable.
If a proper sub SF $\EE'=(E', \Gamdot')\subset \EE=(E,\Gamdot)$
satisfies $p^{\Delta_0}(\EE')=p^{\Delta_0}(\EE)$, then
$\EE'$ is saturated and both $\EE'$ and $\EE/ \EE'$ are
$\Delta_0$-stable.
\end{prop}
\begin{pf}
Remark that if a SF $\EE$ of type $\Hf$ is semistable with respect to
some parameter $\Delta$ then $\EE$ becomes full.
$\EE'$ clearly is saturated and $\Delta_0$-semistable.
If $\EE'$ is not $\Delta_0$-stable, then there is a proper sub SF
$\EE''=(E'', \Gamdot'')$ of $\EE'$ such that
$p^{\Delta_0}(\EE)=p^{\Delta_0}(\EE')=p^{\Delta_0}(\EE'')$.
This implies both $W(\EE',\Hf)$ and $W(\EE'', \Hf)$ are SF-wall
containing $\Delta_0$, so $W(\EE',\Hf)$ equals $W(\EE'', \Hf)$.
Thus we find a constant $\lambda$ such that
\[ p^{\Delta}(\EE)-p^{\Delta}(\EE')= \lambda\left\{ p^{\Delta}(\EE)-
p^{\Delta}(\EE'')  \right\}  \]
for all $\Delta$.
One can deduce that
\[ \frac{i}{\rk E}-\frac{\rk\Gamma'_i}{\rk E'} =
   \lambda \left\{ \frac{i}{\rk E}- \frac{\rk\Gamma''_i}{\rk E''} \right\}  \]
for all $i$, which means that
\begin{equation}\label{eq:alwaysJHF}
  \frac{1}{\rk E}- \frac{\rk( \Gamma'_i / \Gamma'_{i-1})}{\rk E}=
   \lambda \left\{  \frac{1}{\rk E}-
\frac{\rk (\Gamma''_i / \Gamma''_{i-1})}{\rk E''}   \right\} 
\end{equation}
for all $i$.
Since $\EE$ is full, $\rk( \Gamma'_i / \Gamma'_{i-1})$ is either $0$ or $1$.
If $\rk( \Gamma'_i / \Gamma'_{i-1})$ equals $1$ for all $i$ then
it follows that $H^0(E')=H^0(E')\cap\Gamma_n=\Gamma'_n=\Gamma_n=H^0(E)$.
This is contradiction since $E$ is generated by global sections from
\eqref{eq:(O)}. Accordingly 
\begin{equation}\label{eq:somei0}
\rk( \Gamma'_{i_0} / \Gamma'_{i_0-1})=0 \quad\text{for some $i_0$}.
\end{equation}
As to this $i_0$, one can check that
\begin{equation}\label{eq:thisi0}
\rk( \Gamma''_{i_0} / \Gamma''_{i_0-1})=1.
\end{equation}
From \eqref{eq:alwaysJHF}, \eqref{eq:somei0} and \eqref{eq:thisi0}
we can determine $\lambda$ and hence show that
\begin{equation}\label{eq:alli}
\rk(\Gamma'_i/ \Gamma'_{i-1})+ \rk(\Gamma''_i/ \Gamma''_{i-1})=1
\quad\text{for all $i$} .
\end{equation}
On the other hand, $H^0(E'')\neq 0$ by \eqref{eq:(O)}, so there
is a nonzero section $\tau\in H^0(E'')$.
Since $\EE$ is full, some $j$ enjoys the property that
\[ \tau\not\in H^0(E'')\cap \Gamma_{j-1}= \Gamma''_{j-1} \quad\text{and that}
   \quad \tau\in H^0(E'')\cap \Gamma_j =\Gamma''_j. \]
As to this $j$, it also holds that
\[  \tau\not\in H^0(E')\cap \Gamma_{j-1}= \Gamma'_{j-1} \quad\text{and that}
   \quad \tau\in H^0(E')\cap \Gamma_j =\Gamma'_j. \]
However these facts contradict \eqref{eq:alli}. Therefore $\EE'$ is
$\Delta_0$-stable.
\end{pf}
\begin{cor}\label{cor:FlFr}
Assume that $\Hf$ has the property $(A)$,
two parameters $\Delta_-$ and $\Delta_+$ are
contained in adjacent SF-chambers of type $\Hf$, and that
$\Delta_0= t\Delta_- + (1-t)\Delta_+$ $(0<t<1)$ is contained in just
one SF-wall of type $\Hf$.
Suppose that a SF $\EE$ of type $\Hf$ is $\Delta_-$-semistable and
not $\Delta_+$-semistable and hence there is an exact sequence of SFs
\[ 0 \longrightarrow \FFl \longrightarrow \EE \longrightarrow \FFr 
\longrightarrow 0, \]
where $\FFl$ is saturated and satisfies
$p^{\Delta_+}(\FFl)>p^{\Delta_+}(\EE)$.
$($We call such a sub SF $\FFl$ a $\Delta_+$-{\normalshape destabilizer} of
$\EE$.$)$ Then the following holds:
\begin{enumerate}
 \item $\EE$ is $\Delta_-$-stable, and its $\Delta_+$-destabilizer is unique.
 \item If a SF $\EE'$ is endowed with a nontrivial exact sequence
\[ 0 \longrightarrow \FFr \longrightarrow \EE' \longrightarrow \FFl 
\longrightarrow 0, \]
then $\EE'$ is $\Delta_+$-semistable.
\end{enumerate}
\end{cor}
\begin{pf}
The definition of SF-chambers deduces that
$p^{\Delta_0}(\FFl)=p^{\Delta_0}(\EE)$ and that
any $\Delta_+$-destabilizer $\FF$ of $\EE'$, if it exists, satisfies
that $p^{\Delta_0}(\FF)=p^{\Delta_0}(\EE')$.
This corollary follows these facts and the lemma above.
\end{pf}
\section{Moduli theory of SFs}
\label{ss:moduli}
Let us begin with the construction of the coarse moduli scheme of
semistable SFs. 
We fix $\Hf=\left( \Hf'=(f,f_0,f_1), r_1,\dots, r_n\right)$.
\begin{defn}
For a scheme $S$ over $k$, a {\it $S$-flat family of SFs} on $X$ is 
a pair $(E_S, \{ \Gamma_{i,S}\}_i)$ consisting of a $S$-flat sheaf $E_S$
on $X_S$ and a sequence of quotients
\[ Ext_{X_S/S}^2(E_S,K_X) \twoheadrightarrow \left(\Gamma_{n,S}\right)^{\vee}
\twoheadrightarrow \dots \twoheadrightarrow 
\left( \Gamma_{1,S}\right)^{\vee},  \]
where $\Gamma_{i,S}$ is a locally-free ${\cal O}_S$-module.
{\it A homomorphism} $f:\EE'_S=(E'_S,\Gamma'_{\bullet,S})\rightarrow
\EE_S=(E_S,\Gamma_{\bullet, S})$ {\it of flat families of SFs} of length
$n$ is a homomorphism $f:E'_S\rightarrow E_S$ which induces a
homomorphism $f_i:\Gamma'_{i,S}\rightarrow \Gamma_{i,S}$ ($1\leq i\leq n$)
such that
\begin{equation*}
\xymatrix{
Ext^2_{X_S/S}(E_S,K_X) \ar[r]^f \ar[d] & Ext^2_{X_S/S}(E'_S,K_X) \ar[d]
 \\
(\Gamma_{i,S})^{\vee} \ar[r]^{f_i^{\vee}} & (\Gamma'_{i,S})^{\vee}}
\end{equation*}
is commutative. 
\end{defn}
Define a functor $\underline{M}:
(\text{Sch}/k)^{\circ} \rightarrow (\text{Sets})$ as follows.
We first set $\underline{M}'(S)$ is to be the 
set of all $S$-flat families of $\Delta$-semistable SFs of type $\Hf$ 
on $X$ and
$\underline{M}(S)$ is the quotient
$\underline{M}'(S)/\sim$, where $S$-flat families
$\EE_S$ and $\FF_S$ are equivalent if and only if it holds that
$\EE_S\otimes L \simeq \FF_S$ for some line bundle $L$ on $S$.
We also define a functor $\underline{M}^s$ by
replacing ``$\Delta$-semistable'' with ``$\Delta$-stable'' here.
\begin{prop}\label{prop:moduli}
The functor $\underline{M}$ has the coarse moduli
scheme $M(\Delta,\Hf)$ which is projective over $k$.
$M(\Delta, \Hf)(k)$ coincides with the set of all S-equivalence
classes of $\Delta$-semistable SFs of type $\Hf$.
Some open subset $M^s(\Delta,\Hf)\subset M(\Delta, \Hf)$
is the coarse moduli scheme of the functor
$\underline{M}^s$.
\end{prop}
\begin{pf}
One can prove this proposition in a similar fashion to Simpson's 
construction of moduli schemes of semistable sheaves \cite{Sim:moduli}
and \cite[Chap. 4]{HL:text}. We also take the case of parabolic sheaves
\cite{Mar-Yok:parabolic} and of coherent systems \cite{He:CohSyst}
as models.\par
First, there is an integer $m$ such that if
$F$ belongs to ${\cal S}_2(\Hf')\cup {\cal S}_3(\Hf')$, then
both $F$, $F\otimes L$, $F\otimes L(-C)$ and $L$ have the property $(O_m)$.
Let $V_m$ be a $f_1(m)$-dimensional vector space, and denote
by $Q(m,\Hf')$
%$\Quot\left(V_m\otimes L^{-1}(-m), \Hf'\right)=Q(m,\Hf')$ 
Grothendieck's Quot-scheme parametrizing quotient ${\cal O}_X$-modules
of $V_m\otimes L^{-1}(-m)$ whose type is $\Hf'$, and by
$U\subset Q(m,\Hf')$ the open subset of all quotients
$q: V_m\otimes L^{-1}(-m)\twoheadrightarrow E$ such that $E$ is
torsion-free, both $E$, $E\otimes L$ and $E\otimes L(-C)$ have the
property $(O_m)$, and $H^0(q): V_m \rightarrow H^0(E\otimes L(m))$ is
injective.
$Q(m,\Hf')$ has a universal family 
$V_m\otimes {\cal O}_{X_Q} \twoheadrightarrow E_Q\otimes L(m)$ on
$X_{Q(m,\Hf')}$.\par
Next, consider a functor
$\underline{Fl}\left( Ext_{X_U/U}^2(E_U,K_X),r_{\bullet} \right):
(\text{Sch}/U)^{\circ} \rightarrow (\text{Sets})$
which associates with $S\rightarrow U$ the set of all sequences
of surjective homomorphisms
\[ Ext_{X_U/U}^2(E_U,K_X)\underset{U}{\otimes} {\cal O}_S 
\twoheadrightarrow \left(\Gamma_{n,S}\right)^{\vee}
\twoheadrightarrow \dots \twoheadrightarrow 
\left( \Gamma_{1,S}\right)^{\vee} \]
consisting of locally-free ${\cal O}_S$-modules $\Gamma_{i,S}$ with
rank $r_i$.
This is represented by a $U$-scheme, say $R_m$.
By the choice of $U$ a natural map
$Ext^2_{X_U/U}(E_U\otimes L(m),K_X)\otimes H^0(L(m))\rightarrow
Ext^2_{X_U/U}(E_U,K_X)$ is surjective and
$Ext^2_{X_U/U}(E_U\otimes L(m),K_X)\rightarrow
Ext^2_{X_U/U}(V_m\otimes {\cal O}_{X_U},K_X) \simeq V_m^{\vee}\otimes 
{\cal O}_U$ is isomorphic.
Thus if we put $B_m=H^0(L(m))$ then $R_m$ is embedded in
$U\times Fl(V_m^{\vee}\otimes B_m,r_{\bullet})$, where
$Fl(V_m^{\vee}\otimes B_m,r_{\bullet})$ is the flag scheme parametrizing
sequences of surjective maps
$V_m^{\vee}\otimes B_m \twoheadrightarrow \Gamma_n^{\vee} \twoheadrightarrow
\dots \twoheadrightarrow \Gamma_1^{\vee}$ consisting of vector spaces
$\Gamma_i$ with rank $r_i$.\par
Last, a natural map
$Ext^2_{X_U/U}(E_U\otimes L(m)|_C, K_X) \rightarrow Ext^2_{X_U/U}(E_U\otimes
L(m),K_X) \simeq V_m^{\vee}\otimes {\cal O}_U$ induces a morphism
$U \rightarrow Gr\left(V_m,f_1(m)-f_0(m)\right)=:Gr (f_1-f_0)$
to the Grassmannian parametrizing quotient vector spaces of $V_m$ with
rank $f_1(m)-f_0(m)$.
Hence we obtain a embedding
$R_m \subset U\times Fl(V_m^{\vee}\otimes B_m,r_{\bullet}) \times
Gr(f_1-f_0)$ and its closure
\begin{equation}\label{eq:Rbar}
\overline{R_m}\subset Q(m,\Hf') \times Fl(V_m^{\vee}\otimes 
B_m,r_{\bullet}) \times Gr(f_1-f_0)
\end{equation}
which is invariant under the natural action of $G=\SL(V_m)$.\par
Remember some $G$-linearized line bundles on the right side of
\eqref{eq:Rbar};
\[{\cal O}^l_{Q,0}(1)= \det Rp_{X\,*}(E_Q\otimes L(-C)(l))
\text{ and }
{\cal O}^l_{Q,1}(1)= \det Rp_{X\,*}(E_Q\otimes L(l)) \]
are $G$-linearized ample line bundles on $Q(m,\Hf')$ when $l$ is
sufficiently large \cite[Prop. 2.2.5]{HL:text}.
$Fl(V_m^{\vee}\otimes B_n,r_{\bullet})$ has a universal family
\begin{equation}\label{eq:UnivOfFl}
V_m^{\vee}\otimes B_m \otimes {\cal O}_{Fl} \twoheadrightarrow 
\left(\Gamma_{n,Fl}\right)^{\vee} \twoheadrightarrow
\dots \twoheadrightarrow \left(\Gamma_{1,Fl}\right)^{\vee}.
\end{equation}
For positive integers $k_1,\dots,k_n$,
\[{\cal O}_{Fl}(k_1,\dots,k_n)=\sideset{}{_{i=1}^n}\otimes
 \left( \det \Gamma_{i,Fl}^{\vee} \right) ^{\otimes k_i} \]
is a $G$-linearized ample line bundle on
$Fl(V_m^{\vee}\otimes B_m, r_{\bullet})$ by the Pl\"{u}cker embedding.
Similarly, if we denote a universal family of $Gr(f_1-f_0)$ by
$V_m\otimes {\cal O}_{Gr}\twoheadrightarrow W_{Gr}$, then
${\cal O}_{Gr}(1)=\det W_{Gr}$ is a $G$-linearized ample line bundle
on $Gr(f_1-f_0)$.
For a parameter $\Delta=(a, \delta_{\bullet})$ we put
$ f^{\Delta}(l)= (1-a)f_0(l)+af_1(l)+\sum_{i=1}^n \delta_i\cdot r_i $
and then
\begin{multline*}
 L_l= {\cal O}_{Q,0}^l\left( (1-a)f^{\Delta}(m) \right)\otimes
        {\cal O}_{Q,1}^l\left( af^{\Delta}(m) \right) \otimes \\
        {\cal O}_{Fl}\left(\delta_1(f^{\Delta}(l)-f^{\Delta}(m)),\dots,
         \delta_n(f^{\Delta}(l)-f^{\Delta}(m))  \right) \otimes
        {\cal O}_{Gr}((1-a)f^{\Delta}(l)) 
\end{multline*}
is a $G$-linearized $\QQ$-ample line bundle on $R_m$ when $l$ is
sufficiently large.
The GIT quotient $\overline{R_m}^{ss}(L_l) //G$ 
($\overline{R_m}^s(L_l)/G$, resp.) is the
moduli scheme $M(\Delta,\Hf)$ ($M^s(\Delta,\Hf)$, resp.)
if $m$ is sufficiently large
and if $l$ is sufficiently large with respect to $m$.
Its proof proceeds in a similar fashion to that of Theorem 4.3.3
in \cite{HL:text}, so is left to the reader.
\end{pf}
\begin{prop}
If $\Hf$ has the property $(A)$, then $M^s(\Delta,\Hf)$
represents the functor $\underline{M}^s$.
\end{prop}
\begin{pf}
The sheaves $E_U$ and $\Gamma_{i,Fl}$ in the proof of Proposition
\ref{prop:moduli} give a flat family of SFs $\EE_{R^s}$ over
$R^s:=R^s_m(L_l)$.
On the other hand one can check that
$R^s\rightarrow \Ms=M^s(\Delta,\Hf)$
is a $PGL(V_m)$-bundle in a similar fashion to
Proposition 6.4 in \cite{Ma:moduli2}. 
Because $\lambda\cdot\id \in GL(V_m)$ acts on the line bundle
$\Gamma_{1,Fl}$ at \eqref{eq:UnivOfFl} by the multiplication of
$\lambda$,
\begin{multline}
 \EE_{R^s}\otimes \Gamma_{1,Fl}^{\vee} =
 \left( E_U\otimes \Gamma_{1,Fl}^{\vee}\, , \right. \\
 \left.
 Ext^2_{X_{R^s}/R^s}(E_U\otimes \Gamma_{1,Fl}^{\vee}, K_X) \rightarrow
 \Gamma_{n,Fl}^{\vee}\otimes \Gamma_{1,Fl} \rightarrow 
 \dots \rightarrow \Gamma_{1,Fl}^{\vee}\otimes \Gamma_{1,Fl}    \right)
\end{multline}
descends to a $M^s(\Delta,\Hf)$-flat family 
$\overline{\EE}_{\Ms} = (\bE_{\Ms}, \bGam_{\bullet,\Ms})$ from
fpqc descent theory.
By the assumption $p_{X@,*}(\bE_{\Ms})=:\bV_{\Ms}$ is a vector bundle on
$\Ms$ endowed with a flag structure. After \cite[p. 49]{HL:text} we
denote by
$Hom_{\Ms}^-(\bV_{\Ms},\bV_{\Ms}) \subset Hom_{\Ms}(\bV_{\Ms},\bV_{\Ms}) $
the subsheaf consisting of all homomorphisms which preserves the flag
structure, and by $Hom_{\Ms}^+(\bV_{\Ms},\bV_{\Ms})$ its quotient.
By $\id\in \operatorname{End}(E_{\Ms})$ and a natural map
\begin{equation}\label{eq:MapHom+}
  RHom_{X_{\Ms}/\Ms}(\bE_{\Ms},\bE_{\Ms}) \longrightarrow
Hom_{\Ms}(\bV_{\Ms},\bV_{\Ms}) \longrightarrow Hom_{\Ms}^+(\bV_{\Ms},\bV_{\Ms})
\end{equation}
we obtain two triangles 
\begin{equation}\label{eq:TrId}
{\cal O}_{\Ms} \longrightarrow K_{\Ms}^0[-1] \longrightarrow Mc(\id) 
\overset{+1}{\longrightarrow} \quad\text{and}
\end{equation}
\begin{equation}\label{eq:TrHom+}
 RHom_{X_{\Ms}/\Ms}(\bE_{\Ms},\bE_{\Ms}) \longrightarrow
Hom_{\Ms}^+(\bV_{\Ms},\bV_{\Ms}) \longrightarrow K_{\Ms}^0
\overset{+1}{\longrightarrow}.
\end{equation}
\begin{clm}\label{clm:simple}
${\cal H}^0(\id): {\cal O}_{\Ms}\rightarrow {\cal H}^0({\cal K}^0_{\Ms}[-1])
\simeq \HomSF (\overline{\EE}_s,\overline{\EE}_s)$ is isomorphic.
\end{clm}
\begin{pf}
Let $s$ be a point in $\Ms$. From triangles
$\text{\eqref{eq:TrId}}\dotimes k(s)$ 
and $\text{\eqref{eq:TrHom+}}\dotimes k(s)$ one can check that
%\[{\cal H}^{-1}( {\cal K}_{\Ms}^0)\simeq \Ker\left(
%\Hom_X(\bE_s,\bE_s)\rightarrow \Hom^+(\Gamma(\bE_s),\Gamma(\bE_s))\right)\]
%and that
${\cal H}^i\left( Mc(\id)\dotimes k(s) \right)=0$ when $i\leq 0$
since $\overline{\EE}_{\Ms}$ is a flat family of $\Delta$-stable and
accordingly simple SFs.
Because $RHom_{X_{\Ms}/\Ms}(\bE_{\Ms},\bE_{\Ms})$ is perfect, also
${\cal K}_{\Ms}^0$ is. Thus \cite[Thm. 22.5]{Mat:text} verifies
${\cal H}^i(Mc(\id)))=0$ when $i\leq 0$.
\end{pf}
By this claim, this proposition is shown similarly to
\cite[Prop. 4.6.2.]{HL:text}.
\end{pf}
Here we mention the infinitesimal deformation of a SF 
$\GG=(G,\Gamma_{\bullet})$ which satisfies $H^i(G)=0$ when $i>0$;
it is a variation of the standard deformation theory of sheaves
(\cite{Mk:symplectic}, \cite{Ty:symplectic} and others).
Define a functor ${\cal D}$ from the category of Artinian local
$k$-algebras to that of sets by
\[ {\cal D}(A)= \left.
\left\{ \GG_A \bigm| \text{an $A$-flat family of SFs such that
$\GG_A\otimes k\simeq \GG$ } \right\} \right/ \simeq \]
and ${\cal D}(f:A\rightarrow A')(\GG_A)=f^* \GG_A$.
If we put $V=H^0(G)$, we have the following commutative diagram whose 
rows and columns are triangles:
\begin{equation}\label{eq:SixTr}
\xymatrix{
 k \ar[r]^{\id} \ar[d] & RHom_X(G,G) \ar[r] \ar[d]^{\phi} & 
 RHom_X(G,G)/k \ar[r]^-{+1} \ar[d]^{\phi_+} & \\
 Hom^-(V,V) \ar[r] \ar[d] & Hom(V,V) \ar[r] \ar[d] &
 Hom^+(V,V) \ar[r]^-{+1} \ar[d] & \\
 Hom^-(V,V)/k \ar[d]^{+1} \ar[r] & Mc(\phi) \ar[r]^{\alpha} \ar[d]^{+1} &
Mc(\phi_+) \ar^{+1}[r], \ar[d]^{+1} & \\
 & & & }
\end{equation}
where $\phi$ is the map \eqref{eq:MapHom+}.
Let $A\rightarrow \bA$ be a small extension of Artinian 
local rings, that is, a surjective ring homomorphism whose kernel
${\frak a}$ satisfies ${\frak a}\cdot m_A=0$.
\begin{lem}\label{lem:ob}
Let $\GG_A$ and $\GG'_A$ be elements in ${\cal D}(A)$ endowed with an
isomorphism
$\overline{\kappa}: \GG_A\otimes\bA \simeq \GG'_A\otimes\bA$. 
Then there is an obstruction class
$\ob(\overline{\kappa},{\frak a})\in {\cal H}^0(Mc(\phi_+))\otimes {\frak a}$
with the property that $\ob=0$ if and only if
$\overline{\kappa}$ extends to an isomorphism
$\kappa: \GG_A \simeq \GG'_A$.
Conversely, let $\GG_A$ be an $A$-flat family of SFs extending $\GG$.
For any $v\in H^0(Mc(\phi_+))\otimes {\frak a}$ we have an $A$-flat
family of SFs $\GG'_A$ and an isomorphism
$\overline{\kappa}:\GG_A\otimes\bA \simeq \GG'_A\otimes\bA$
such that $\ob(\overline{\kappa},{\frak a})=v$.
\end{lem}
\begin{pf}
We shall utilize methods in \cite{Lau:Massey} or
\cite[Section 2.A.6]{HL:text}.
The sheaf $G_0:= G'_A\otimes k$ has an injective resolution
$ 0 \rightarrow G_0 \overset{\epsilon_0}{\rightarrow} I^0
 \overset{d_0}{\rightarrow} I^1 \rightarrow \dots$.
One can find an exact sequence
\begin{equation*}
0 \longrightarrow G'_A \overset{\epsilon'}{\longrightarrow} I^0\otimes A 
 \overset{d'_A}{\longrightarrow} I^1\otimes A \longrightarrow \dots
\end{equation*}
such that $\epsilon'\otimes k=\epsilon_0$ and $d'_A\otimes k=d_0$,
and an exact sequence
\begin{equation*}
0 \longrightarrow G_A \overset{\epsilon}{\longrightarrow} I^0\otimes A 
 \overset{d_A}{\longrightarrow} I^1\otimes A \longrightarrow \dots
\end{equation*}
such that $(\epsilon'\otimes\bA)\circ\overline{\kappa}=\epsilon\otimes \bA$
and $d_A\otimes \bA= d'_A\otimes \bA$. Then
$\partial=d_A-d'_A:I^{\bullet} \rightarrow I^{\bullet +1}\otimes {\frak a}$
lies in
$Z^1\left( \Hom_X^{\bullet}(I^{\bullet},I^{\bullet})\otimes {\frak a}\right)$.
Since ${\cal H}^1\left( \Hom^{\bullet}(\Gamma(I^{\bullet}),\Gamma(I^{\bullet}))
\otimes {\frak a}\right)= \Ext^1(V,V)\otimes{\frak a}$ is zero,
$\Gamma(\partial)$ belongs to $B^1(\Hom^{\bullet}(\Gamma(I^{\bullet}),
\Gamma(I^{\bullet}))\otimes{\frak a})$, in other words,
$\Gamma(\partial)=-\Gamma(d)e+e\Gamma(d)$ with some
$e\in\Hom^0(\Gamma(I^{\bullet}), \Gamma(I^{\bullet}))\otimes{\frak a}$.
Therefore, the diagram
\begin{equation*}
\xymatrix{ 
\Gamma(I^{\bullet}\otimes A) \ar[r]^{\Gamma(d_A)} \ar[d]_{1-e} &
\Gamma(I^{\bullet +1}\otimes A) \ar[d]^{1-e} \\
\Gamma(I^{\bullet}\otimes A) \ar[r]^{\Gamma(d'_A)} & 
\Gamma(I^{\bullet +1}\otimes A)}
\end{equation*}
is commutative and induces a map 
$1-e: p_{X*}(G_A)\rightarrow P_{X*}(G'_A)$.
We can choose $e$ so that this $1-e$ commutes with flag structures
because $\Gamma(\overline{\kappa})$ does.
One can verify that
\[ (-e,\partial)\in Z^0\left(Mc(p_{X*}: \Hom^{\bullet}_X(I^{\bullet}
,I^{\bullet}) \rightarrow 
\Hom^{\bullet}(\Gamma(I^{\bullet}),\Gamma(I^{\bullet})))\right)
\otimes{\frak a} \]
and hence obtains
$[(-e,\partial)]\in {\cal H}^0(Mc(p_{X*}))\otimes{\frak a}\simeq 
{\cal H}^0(Mc(\phi))\otimes {\frak a}.$
Its image by $\alpha$ in \eqref{eq:SixTr},
$\alpha[(-e,\partial)]\in {\cal H}^0(Mc(\phi_+))\otimes{\frak a}$,
is independent of the choice of $d_A$, $d'_A$ and $e$, and
equals zero if and only if $\overline{\kappa}$ extends to an isomorphism
$\kappa: \GG_A \simeq \GG'_A$; its proof is left to the reader.
As to the ``Conversely'' part, one can prove it by reversing the
construction above.
\end{pf}
\begin{cor}\label{cor:TgtMm}
Let $\pi:A\rightarrow\bA$ be a small extension.
\begin{enumerate}
 \item For $\GG_{\bA}\in {\cal D}(\bA)$, there is a class
$\ob(\GG_{\bA},{\frak a})\in {\cal H}^1(Mc(\phi_+))\otimes{\frak a}
\simeq \Ext^2_X(G,G)\otimes{\frak a}$ with the property that
$\ob=0$ if and only if some $\GG_A\in{\cal D}(A)$ satisfies
$\GG_A\otimes\bA\simeq \GG_{\bA}$.
 \item Suppose $\GG$ is a simple SF. Then the fiber
${\cal D}(\pi)^{-1}(\GG_{\bA})$ is an affine space with the
transformation group ${\cal H}^0(Mc(\phi_+))\otimes{\frak a}$ unless
it is empty.
\end{enumerate}
 \end{cor}
\begin{pf}
 Since $p_{X*}(G_{\bA})$ is a locally-free $\bA$-module, (i) follows
from deformation theory of sheaves. (ii) results from Claim
\ref{clm:simple} and Lemma \ref{lem:ob}.
\end{pf}
\section{Set of $\Delta_-$-semistable and not $\Delta_+$-semistable SFs}
\label{ss:SetP}
We hereafter assume that $\Hf$ has the property $(A)$, and parameters
$\Delta_{\pm}$ and $\Delta_0$ meet the conditions in Corollary
\ref{cor:FlFr}.
$M_+$ and $\Mm$ mean $M(\Delta_+,\Hf)$ and
$M(\Delta_-,\Hf)$ for short.
Since $M(\Delta_-,\Hf)=M^s(\Delta_-,\Hf)$ by the remark
after Definition \ref{defn:(A)}, there is a functor
$\underline{P}: (\operatorname{Sch}/\Mm)^{\circ}\rightarrow 
(\operatorname{Sets})$ 
which associates $q:S\rightarrow \Mm$ with the set of
all isomorphic classes of $S$-flat families
$\tau:\FFl_S\rightarrow q^* \overline{\EE}_{\Mm}$ of
$\Delta_+$-destabilizers, that is,
$\tau$ is a homomorphism of flat families of SFs such that, for any
point $s\in S$, $\tau\otimes k(s): \FFl_s\rightarrow \overline{\EE}_s$
gives a $\Delta_+$-destabilizer of $\overline{\EE}_s$.
\begin{lem}
A closed subscheme $P\subset\Mm$ represents the functor $\underline{P}$.
\end{lem}
\begin{pf}
$\underline{P}$ is representable by Grothendieck's Quot-schemes and
\cite[Lem. 3.1]{Mar-Yok:parabolic}.
If $\FFl$ is a $\Delta_+$-destabilizer of a $\Delta_-$-semistable
SF $\EE$ of type $\Hf$, then $\HomSF (\FFl,\EE/\FFl)=0$ by Proposition
\ref{prop:Rk2-like}. Hence the same argument as in the proof of
\cite[Lem. 2.2]{Yam:Dthesis} shows this lemma.
\end{pf}
There are $P$-flat families of SFs 
$\FFl_P=(\Fl_P, \Gaml_{\bullet,P})$ and
$\FFr_P=(\Fr_P, \Gamr_{\bullet,P})$, 
and an exact sequence of families of SFs
\begin{equation}\label{eq:HNFforSF}
 0 \longrightarrow \FFl_P \longrightarrow
\overline{{\cal E}}_{\Mm}|_P \longrightarrow \FFr_P \longrightarrow 0. 
\end{equation}
Now let $T\subset \overline{T}$ be a closed immersion whose ideal sheaf
${\frak a}\subset{\cal O}_{\Mm}$ satisfies that ${\frak a}^2=0$, and
$f: \overline{T}\rightarrow \Mm$ a morphism such that 
its restriction to $T$ factors through $P\subset \Mm$, in other words,
$f|_T$ induces a morphism $g:T\rightarrow P$.
When we denote $f^*(\overline{E}_{\Mm})=\bE_{\bT}$,
$g^* \Fl_P=\Fl_T$ and so on, the exact sequence of ${\cal O}_{X_{\bT}}$-modules
associated with \eqref{eq:HNFforSF} gives a diagram in $\Coh(X_{\bT})$
\begin{equation}
\xymatrix{
  & & 0 \ar[d] & & \\
 0 \ar[r] & {\frak a}\otimes_T \Fl_T \ar[r] & {\frak a}\otimes_T 
\bE_{\bT} \ar[r] \ar[d] & {\frak a}\otimes_T \Fr_T
 \ar[r] & 0 \\
  & & \bE_{\bT} \ar[d] & & \\
 0 \ar[r] &  \Fl_T \ar[r] & \bE_{\bT}|_T \ar[r] \ar[d]
 &  \Fr_T \ar[r] & 0 \\
  & & 0 & & }
\end{equation}
whose rows and columns are exact. This diagram induces the following: \\
(i) An ${\cal O}_{X_T}$-module $W_T=\Ker( \bE_{\bT} \rightarrow \Fr_T)/ 
\IIm({\frak a}\otimes \Fl_T \rightarrow \bE_{\bT})$ 
and an exact sequence
 \begin{equation}\label{eq:WT}
 0 \longrightarrow {\frak a}\otimes \Fr_T \longrightarrow W_T
 \longrightarrow \Fl_T \longrightarrow 0.
 \end{equation}
Similarly, homomorphisms of ${\cal O}_T$-modules with flag structures
associated with \eqref{eq:HNFforSF} brings following elements:\\
(ii) An ${\cal O}_T$-module 
$\Lambda_{\bullet,T}=\Ker(  \overline{\Gamma}_{\bullet ,\bT} \rightarrow 
\Gamr_{\bullet,T}) / \IIm({\frak a}\otimes \Gaml_{\bullet, T}\rightarrow 
\overline{\Gamma}_{\bullet ,\bT})$ 
and an exact sequence
 \begin{equation}\label{eq:LambdaT}
 0 \longrightarrow {\frak a}\otimes \Gamr_{\bullet, T}
 \longrightarrow \Lambda_{\bullet,T}
 \longrightarrow \Gaml_{\bullet, T} \longrightarrow 0\,;
 \end{equation}
(iii)Homomorphisms $\tau_{\bullet -1}: \Lambda_{\bullet -1,T}\rightarrow
\Lambda_{\bullet, T}$ 
and $\iota_{\bullet}: \Lambda_{\bullet ,T}\rightarrow p_{X*}(W_T)$ such
that the diagram
\begin{equation*}
\xymatrix{
 {\frak a}\otimes \Gamr_{\bullet -1,T} \ar[d] \ar@{^{(}-}[r] &
 {\frak a}\otimes \Gamr_{\bullet ,T} \ar[d] \ar@{^{(}-}[r] &
 {\frak a}\otimes p_{X*}(\Fr_T) \ar[d] \\
 \Lambda_{\bullet -1,T} \ar[r]^{\tau_{\bullet -1}} \ar[d] &
 \Lambda_{\bullet ,T} \ar[r]^{\iota_{\bullet}} \ar[d] &
 p_{X*}(W_T) \ar[d] \\
 \Gaml_{\bullet -1,T} \ar@{^{(}-}[r]  &
 \Gaml_{\bullet, T} \ar@{^{(}-}[r] &
 p_{X*} (\Fl_T), }
\end{equation*}
which is a combination of $p_{X*}(\text{\eqref{eq:WT}})$,
\eqref{eq:LambdaT} and flag structures of $\FFl_P$ and $\FFr_P$, is
commutative.
\begin{lem}
$f:\overline{T}\rightarrow \Mm$ factors through $P\subset\Mm$ if and
only if there are a section $\kappa: \Fl_T \rightarrow W_T$ of
\eqref{eq:WT} and a section
$\kappa_{\bullet}:\Gaml_{\bullet,T}\rightarrow \Lambda_{\bullet ,T}$
of \eqref{eq:LambdaT} which make the following diagram commutative.
\begin{equation*}
\xymatrix{
 \Lambda_{\bullet -1,T} \ar[r]^{\tau_{\bullet -1}} &
 \Lambda_{\bullet ,T} \ar[r]^{\iota_{\bullet}} &
 p_{X*}(W_T) \\
 \Gaml_{\bullet -1,T} \ar@{^{(}-}[r]  \ar[u]_{\kappa_{\bullet -1}} &
 \Gaml_{\bullet, T} \ar@{^{(}-}[r] \ar[u]_{\kappa_{\bullet}}  &
 p_{X*} (\Fl_T) \ar[u]_{p_{X*}(\kappa)}}
\end{equation*}
\end{lem} 
\begin{pf}
It is a variation of the deformation theory of Quot-schemes
\cite[page 43]{HL:text}, so we omit the proof.
\end{pf}
We shall denote $\Vl_P=p_{X*}(\Fl_P)$, $\Vr=p_{X*}(\Fr_P)$ and $\omega_X=K_X[2]$.
By the duality theorem \cite{Ha:residue}, a natural morphism
\begin{equation}\label{eq:psi-}
\psi_-:\,  R Hom_{X_P/P}(\Fl_P,\Fr_P) \longrightarrow
 Hom^+_P(\Vl_P,\Vr_P) 
\end{equation}
induces a triangle in $D^b(P)$
\begin{equation}\label{eq:H^0}
 Hom^+_P(\Vl_P,\Vr_P)^{\vee} \longrightarrow Hom_{X_P/P}(\Fr_P,
  \Fl_P(\omega_X)) \longrightarrow {\cal K} \overset{+1}{\longrightarrow}.
\end{equation}
\begin{lem}
There is an obstruction class 
$\ob(f,g)\in \Ext^1_T(Lg^*({\cal K}),{\frak a})$ with the property
 that
$f:\bT\rightarrow\Mm$ factors through $P\subset\Mm$ if and only if $\ob=0$.
\end{lem}
\begin{pf}
The functor $R\Hom_T(Lg^*(?),{\frak a})$ takes \eqref{eq:H^0} to a
 triangle
\begin{multline}\label{eq:dual-triangle}
 R\Hom_T(Lg^*{\cal K},{\frak a}) \longrightarrow
 R\Hom_{X_T}(\Fl_T,\Fr_T \underset{T}{\otimes} {\frak a}) \\
 \overset{\psi_-}{\longrightarrow} R\Gamma_T(Hom_T^+(\Vl_T,
 \Vr_T\otimes{\frak a})) \overset{+1}{\longrightarrow} .
\end{multline}
Let $\epsilon^{(r)}:\Fr_T\otimes{\frak a}\rightarrow (\Ibr,d^{(r)})$ be 
an injective resolution in $\Mod(X_T)$.
As for a ${\cal O}_T$-module $\Vr_T\otimes{\frak a}$ with the filtration
$\Gamma_{1,T}^{(r)}\otimes{\frak a}\subset\dots
 \Gamma_{n,T}^{(r)}\otimes{\frak a} \subset
 \Gamma_{n+1,T}^{(r)}\otimes{\frak a} =V^{(r)}_T\otimes {\frak a}$,
pick an injective resolution
$\operatorname{gr}^i(V^{(r)}_T\otimes{\frak a})\rightarrow (K_i^{\bullet},d_i)$
for $i=1,\dots,n+1$ and find an injective resolution
$\Vr_T\otimes{\frak a} \rightarrow
 (K^{\bullet}=\sideset{}{_{j=1}^{n+1}}\oplus K^{\bullet}_j,d_K)$ such that
$d_K(\sideset{}{_{j\leq i}}\oplus K_j^{\bullet})\subset \sideset{}{_{j\leq
 i}}\oplus K_j^{\bullet}$ 
and that
$\operatorname{gr}^i(d_K):K_i^{\bullet}\rightarrow K_i^{\bullet +1}$ coincides with
 $d_i$ for every $i$.
In particular $(K^{\bullet},d_K)$ is a filtered complex.
One can describe $R\Hom_T(Lg^*{\cal K},{\frak a})$ by $\Ibr$ and $K^{\bullet}$.
Indeed, a natural map
$\Vr_T\otimes{\frak a}\rightarrow p_{X*}(\Fr_T\otimes{\frak a})$ is isomorphic,
and its inverse map extends to a quasi-isomorphism
\begin{equation}\label{eq:nu}
  \nu : (p_{X*}(I^{\bullet}),p_{X*}(d_I)) \longrightarrow (K^{\bullet},d_K).
\end{equation}
Fix an affine open covering $\{ T_a \}$ of $T$ such that the exact
sequence $p_{X*}(\text{\eqref{eq:WT}})|_{T_a}$,
\[ 0 \longrightarrow {\frak a}\otimes \Vr_T|_{T_a} \longrightarrow
   p_{X*}(W_T)|_{T_a} \longrightarrow \Vl_T|_{T_a} \longrightarrow 0, \]
has a section $j_a: \Vl_T|_{T_a}\rightarrow p_{X*F}(W_T)|_{T_a}$ which
preserves filtrations $\Gaml_{\bullet T}$ and $\Lambda_{\bullet T}$.
Since $K^{\bullet}$ has a filtration, we obtain complexes 
$Hom_T^+(\Vl_T,K^{\bullet})$ and
\[(C^{\bullet}(\{T_a\}, Hom_T^+(\Vl_T,K^{\bullet})), 
(-1)^{\operatorname{deg}}d_{Cech}+d_K), \] 
where $(C^{\bullet}(\{T_a\}, Hom_T^+(\Vl_T, K^q)), d_{Cech})$ is the C\v{e}ch
complex. The homomorphism $\nu$ \eqref{eq:nu} derives
\[ p_{X*}(\nu):\, \Hom_{X_T}(\Fl_T,\Ibr) \longrightarrow 
 C^{\bullet}(\{T_a\}, Hom_T^+(\Vl_T,K^{\bullet})), \]
and $R\Hom_T(Lg^*{\cal K},{\frak a})[1]$ at \eqref{eq:dual-triangle} is
represented by $Mc(p_{X*}(\nu))$. \par
Let $\alpha\in Z^1(\Hom_{X_T}(\Fl_T,I^{\bullet}))$ represent the image
of identity map by the map 
$\Hom_{X_T}(\Fl_T,\Fl_T)\rightarrow \Ext^1_{X_T}(\Fl_T,\Fr_T\otimes{\frak a})$
coming from \eqref{eq:WT}.
By the exact sequence \eqref{eq:WT}, 
$\epsilon^{(r)}:{\frak a}\otimes \Fr_T \hookrightarrow I^{0\,(r)}$
extends to $\epsilon': W_T\rightarrow I^{0\,(r)}$.
If we denote $\Vl_T|_{Ta} \overset{j_a}{\rightarrow} p_{X*}(W_T)|_{T_a}
\overset{\epsilon'}{\rightarrow} p_{X*}(I^{0(r)})
 \overset{\nu}{\rightarrow} K^0$ 
by $i_a$, then one can check that 
$(\alpha, \{\overline{i_a}\}) \in [Mc(p_{X*}(\nu))]_0$ is contained in
$Z^0(Mc(p_{X*F}(\nu)))$ and that
\[\ob(f,g):= [(\alpha, \{\overline{i_a}\})]\in
 H^0(Mc(p_{X*}(\nu)))\simeq \Ext^1_T(Lg^*{\cal K},{\frak a}) \]
enjoys the property asserted in this lemma.
\end{pf}
When sheaves $G$ and $G'$ on a scheme $S$ have filtrations
$\{G_i\subset G\}$ and $\{ G'_i\subset G'\}$ of length $n$, we have
objects $R\Hom_S^-(G,G')$ and $R\Hom_X^+(G,G')$ in $D(S)$ with a
triangle 
\begin{equation}\label{eq:Hom+-}
R\Hom_S^-(G,G') \longrightarrow R\Hom_S(G,G') \longrightarrow
R\Hom_S^+(G,G') \overset{+1}{\longrightarrow};
\end{equation}
see \cite[p. 49]{HL:text}.
For a point $s\in P$ corresponding to a SF $\EE$, we here explain how to
derive the following diagram whose rows and columns are triangles:
\begin{equation}\label{eq:Hom+-forSF}
\xymatrix{
R\Hom^{(-)}_X(E,E)/k \ar[r] \ar[d]^{\psi_-} & R\Hom_X(E,E)/k \ar[r]
\ar[d]^{\phi_+} & R\Hom_X(F^{(l)},F^{(r)}) \ar[r]^(.7){+1} \ar[d]^{\psi_+} 
& \\
\Hom^{+(-)}(V,V) \ar[d] \ar[r] & \Hom^+(V,V) \ar[r] \ar[d] &
\Hom^+(V^{(l)},V^{(r)}) \ar[r]^(.7){+1} \ar[d] & \\
Mc(\psi_-) \ar[r] \ar[d]^{+1} & Mc(\phi_+) \ar[r] \ar[d]^{+1} &
Mc(\psi_+) \ar[r]^{+1} \ar[d]^{+1} & \\
 & & &}
\end{equation}
Equation \eqref{eq:HNFforSF} gives filtrations $F^{(l)}\subset E$ and
$\Gamma(F^{(l)})\subset V=\Gamma(E)$, and the flag structures of SFs
are nothing but
filtrations $\Gamma_{\bullet}\subset V$, $\Gamma_{\bullet}^{(l)}\subset
V^{(l)}=\Gamma(F^{(l)})$, and so on.
$R\Hom_X^{(-)}(E,E)$ means $R\Hom^-_X$ with respect to the former
filtration, and the first row in \eqref{eq:Hom+-forSF} comes from \eqref{eq:Hom+-}.
$\Hom^+(V,V)$ means, by the definition, $\Hom^+$ with respect to the
latter filtration, and $\Hom^{+(-)}$ the kernel of a natural map
$\Hom^+(V,V)\rightarrow \Hom^+(V^{(l)},V^{(r)})$.
Morphisms $\phi_+$ and $\psi_-$ are those of \eqref{eq:SixTr} and
\eqref{eq:psi-}, and $\psi_+$ the induced one.
\begin{prop}
The tangent space $T_s P$ is isomorphic to $H^0(Mc(\psi_+))$.
\end{prop}
\begin{pf}
 Since $H^{-1}(Mc(\psi_+))\simeq\HomSF(\FFl,\FFr)=0$, \eqref{eq:Hom+-forSF}
induces an exact sequence
\[ 0 \longrightarrow H^0(Mc(\psi_-)) \longrightarrow H^0(Mc(\phi_+))
 \overset{\varphi}{\longrightarrow} H^0(Mc(\psi_+)). \]
If $f:\bT=\Spec(k[\epsilon]/(\epsilon^2))\rightarrow \Mm$ and
$f':\bT\rightarrow P\subset \Mm$ extend $g=s:T=\Spec k\rightarrow\Mm$,
then $\varphi$ sends
$\ob(\overline{\kappa},k\cdot\epsilon)\in H^0(Mc(\psi_+))\otimes k\cdot\epsilon$
associated with
$\overline{\kappa}:f^*\EE_{\Mm}\otimes k(s)=\EE=f^{'*}\EE_{\Mm}\otimes
 k(s)$ to
$\ob(f,g)\in\Ext^1(Lg^*({\cal K}),k\cdot\epsilon)\simeq
 H^0(Mc(\psi_+))\otimes k\cdot\epsilon$,
where the last equality holds from \eqref{eq:dual-triangle}.
It immediately leads to this proposition.
\end{pf}
\begin{cor}\label{cor:codimPMm}
Let $r$ be an integer and $C$ a compact subset in the ample cone of $X$.
If ${\cal O}(1)$ lies in $C$ and if $s\in P$ corresponds to a SF
$\EE=(E,\Gamma_{\bullet})$ with $\rk(E)=r$, then it holds that
$\operatorname{codim}_s(P,\Mm)\geq \Delta(E)/2r -B(r,X,C)$, 
where $\Delta(E)=2rc_2(E)-(r-1)c_1(E)^2$ and $B(r,X,C)$ is a constant
depending only on $(r,X,C)$.
\end{cor}
\begin{pf}
By the proposition above and Corollary \ref{cor:TgtMm},
\begin{align*}
\dim_s P\leq & \dim H^0(Mc(\psi_+)) \leq 
  \dim\Ext_X^{1(-)}(E,E)-1+ \dim\Hom^+(V,V)\quad\text{and} \\
\dim_s\Mm \geq & \dim H^0(Mc(\phi_+)) - \dim H^1(Mc(\phi_+)) 
 = \dim\Hom^+(V,V) -\chi(E,E)+1.
\end{align*}
Then this corollary results from O'Grady's estimation of
$\dim\Ext^{1(-)}$ (\cite[Prop 3.A.2]{HL:text} and
\cite{OG:basicresults}) and the Riemann-Roch formula.
\end{pf}
\section{Blowing-up construction}\label{ss:blowup}
As Corollary \ref{cor:codimPMm} shows, it is reasonable to expect that
$\Mm$ and $\Mp$ are birationally equivalent.
We here describe how to connect them by a single blowing-up and down in
a moduli-theoretic way.
Let $p:\tilM\rightarrow\Mm$ be the blowing-up along $P$ with exceptional
divisor $E$. Then we have a flat family of SFs
$p^* \bEE_{\Mm}=\bEE_{\tilM}$ over $\tilM$ and an exact sequence of flat
families of SFs
\[ 0 \longrightarrow p^*\FFl_P=\FFl_E \longrightarrow \bEE_{\Mm}|_E
 \longrightarrow \FFr_E \longrightarrow 0 \]
coming from \eqref{eq:HNFforSF}, and we can show the following facts in
the same way as the case of rank-two sheaves (\cite[Section 3 and 4]{Yam:Dthesis})
except for obvious modifications: \\
%
%\begin{enumerate}
(i)$\,\EE'_{\tilM}:=\Ker(\bEE_{\tilM}\rightarrow \bEE_{\tilM}|_E
       \rightarrow \FFr_E)$ is a flat family of SFs over $\tilM$
equipped with an exact sequence
\begin{equation}\label{eq:exactEE'}
  0 \longrightarrow \FFr_E\otimes{\cal O}_E(-E) \overset{k_1}{\longrightarrow}
\EE'_{\tilM}|_E \longrightarrow \FFl_E \longrightarrow 0 
\end{equation}
of families of SFs over $E$. One can regard this as an elementary
transform of families of SFs.\\
(ii) When one applies results in the last section to case where
$f:\Spec({\cal O}_{\tilM}/{\cal O}(-2E))=\bT\overset{p}{\rightarrow}\Mm$
and $g:E=T\overset{p}{\rightarrow}P$, he obtains
${\cal W}_E=(W_E,\Lambda_{\bullet E})$, which is a flat family of SFs since
${\frak a}={\cal O}_E(-E)$ is a line bundle, and an exact sequence
\begin{equation}\label{eq:exactWE}
  0 \longrightarrow \FFr_E\otimes{\cal O}_E(-E) \overset{k_2}{\longrightarrow}
{\cal W}_E \longrightarrow \FFl_E \longrightarrow 0. 
\end{equation}
In fact, there is an isomorphism $\lambda: \EE'_{\tilM}|_E \simeq {\cal W}_E$
of families of SFs which satisfies $\lambda\circ k_1=k_2$ in \eqref{eq:exactEE'}
and \eqref{eq:exactWE}.\\
(iii) For any point $s\in P$, the exact sequence
       $\text{\eqref{eq:exactEE'}}\otimes k(s)$ of SFs is
nontrivial. Consequently $\EE'_{\Mm}$ is a family of
$\Delta_+$-semistable SFs by Corollary \ref{cor:FlFr} and so results in
       a morphism $q:\tilM\rightarrow\Mp$.
%\end{enumerate}
%
\begin{prop}\label{prop:blup}
By reversing $\Delta_-$ and $\Delta_+$ we get a closed subscheme
$P'\subset\Mp$ provided with a similar property to $P\subset\Mm$, and
then the morphism $q:\tilM\rightarrow\Mp$ defined above is the
 blowing-up of $\Mp$ along $P'$. Consequently
\[ \Mm \overset{p}{\longleftarrow} \tilM \overset{q}{\longrightarrow} \Mp \]
are blowing-ups derived from moduli theory.
\end{prop}
We shall end this article with relating variation of parameters $\Delta$
and the $\Delta$-stability of SFs to that of polarizations $H$ on $X$ and
the $H$-semistability of sheaves.
When a class $\cc=(r,c_1,c_2)\in\ZZ\times\NS(X)\times\ZZ$ is given, 
let $H_-$ and $H_+$
be polarizations on $X$ contained in adjacent chambers of type $\cc$,
and $H_0=tH_-+(1-t)H_+$ ($0<t<1$) lie in just one wall of type $\cc$;
see \cite[Def. 2.1]{Yos:chamber} for chambers and walls of type $\cc$.
For positive integers $m,n$ and a constant $0<a<1$,
\[\chi^a(E)(l)=(1-a)\chi(E\otimes mH_-(l))+a\chi(E\otimes nH_+(l))\]
defines {\it $a$-$($semi$)$stability} of a sheaf $E$ on $X$.
\begin{prop}
If $m\gg 0$ and if $n\gg 0$ with respect to $m$, then the following
holds about a sheaf $E$ of type $\cc$:
$E$ is $H_-$-stable if and only if $E$ is $0$-stable if and only
if $E$ is $a_-$-stable where $a_->0$ is sufficiently small.
This also holds when one replaces ``stable'' with ``semistable'' here.
\end{prop}
\begin{pf}
 As to the first ``if and only if'' part, refer to
 \cite[Lem. 3.1]{EG:variation} in rank-two case and \cite[Lem. 3.6]{MW:Mumford}
 for general case. The second is an easy exercise.
\end{pf}
Set ${\cal O}(1)$, $L$ and $C$ to be $H_0$, $nH_+$ and $nH_+-mH_-$
respectively. We can assume that $\Hf=(\chi(E(l)),\chi(E\otimes
L(-C)(l)),\chi(E\otimes L(l)),l_1,\dots,l_n)$, where $E$ is of type
$\cc$, has the property $(A)$.
Choose a parameter $\Delta_{H_-}$ (resp. $\Delta_{H_+}$) so that no
SF-wall of type $\Hf$ separates $\Delta_{H_-}$ and $(0,0,\dots,0)$
(resp. $\Delta_{H_+}$ and $(1,0,\dots,0)$). Then a sheaf $E$ of type
$\cc$ and a SF $\EE=(E,\Gamma_{\bullet})$ of type $\Hf$ satisfies
that (i) $E$ is $H_-$-semistable if $\EE$ is $\Delta_{H_-}$-semistable
and that (ii) $\EE$ is $\Delta_{H_-}$-stable if $E$ is $H_-$-stable.
Thus, if one denotes by
$M(\Delta_{H_-},\cc)\subset M(\Delta_{H_-},\Hf)$ the union
of connected components consisting of all SFs $\EE=(E,\Gamma_{\bullet})$
such that $E$ is of type $\cc$, then one gets a natural morphism to the
coarse moduli scheme $M(H_-,\cc)$ of $H_-$-semistable sheaves of
type $\cc$, $h:M(\Delta_{H_-},\cc)\rightarrow M(H_-,\cc)$,
whose restriction $h:h^{-1}(M^s(H_-,\cc))\rightarrow  M^s(H_-,\cc)$ is 
a Grassmannian-bundle in \'{e}tale topology.
Because there is a sequence of parameters
$\Delta_{H_-}=\Delta_0,\Delta_1,\dots,\Delta_m=\Delta_{H_+}$ such that
$\Delta_i$ and $\Delta_{i+1}$ are in adjacent chambers of type $\Hf$
for all $i$, we arrive at a diagram 
\begin{equation*}
\xymatrix{
 & \tilM_0 \ar[dl]_{p_0} \ar[dr]^{q_0} & &  \dots \ar[dl]_{p_1} &
 \tilM_{m-1} \ar[dr] ^{q_{m-1}} & \\
 M(\Delta_0,\cc) \ar[d]^h & & M(\Delta_1,\cc)  &  & &
 M(\Delta_m,\cc) \ar[d]^h  \\
 M(H_-,\cc) &  & & & & M(H_+,\cc),}
\end{equation*}
where $p_i$ and $q_i$ are blowing-ups in Proposition \ref{prop:blup}.
%
%\bibliography{mybib.bib}
\providecommand{\bysame}{\leavevmode\hbox to3em{\hrulefill}\thinspace}
\providecommand{\MR}{\relax\ifhmode\unskip\space\fi MR }
% \MRhref is called by the amsart/book/proc definition of \MR.
\providecommand{\MRhref}[2]{%
  \href{http://www.ams.org/mathscinet-getitem?mr=#1}{#2}
}
\providecommand{\href}[2]{#2}

\end{document}